\documentclass[preprint,12pt]{elsarticle}
\usepackage{amssymb,amsmath,amsthm}

\biboptions{comma,square,compress}

\journal{Indagationes Mathematicae}

\newcommand*{\mR}{\mathbb R} \newcommand*{\mT}{\mathbb T}
\newcommand*{\mZ}{\mathbb Z}

\newcommand*{\cJ}{\mathcal J} \newcommand*{\cK}{\mathcal K}
\newcommand*{\cL}{\mathcal L} \newcommand*{\cM}{\mathcal M}
\newcommand*{\cO}{\mathcal O} \newcommand*{\cP}{\mathcal P}
\newcommand*{\cT}{\mathcal T} \newcommand*{\cZ}{\mathcal Z}

\newcommand*{\Fi}{\mathfrak i} \newcommand*{\Fj}{\mathfrak j}
\newcommand*{\Fl}{\mathfrak l} \newcommand*{\Fm}{\mathfrak m}
\newcommand*{\Fq}{\mathfrak q} \newcommand*{\Fr}{\mathfrak r}

\newcommand*{\fD}{\mathfrak D} \newcommand*{\fH}{\mathfrak H}
\newcommand*{\fI}{\mathfrak I} \newcommand*{\fK}{\mathfrak K}
\newcommand*{\fL}{\mathfrak L} \newcommand*{\fM}{\mathfrak M}

\newcommand*{\rt}{\mathrm t}

\newcommand*{\vare}{\varepsilon} \newcommand*{\varp}{\varphi}

\newcommand*{\const}{\mathrm{const}}

\newcommand*{\codim}{\mathop{\mathrm{codim}}\nolimits}
\newcommand*{\Fix}{\mathop{\mathrm{Fix}}\nolimits}

\theoremstyle{plain}
\newtheorem{thm}{Theorem}

\theoremstyle{definition}
\newtheorem{dfn}{Definition}
\newtheorem{rem}{Remark}

\begin{document}

\begin{frontmatter}

\title{Hamiltonian and reversible systems \\ with smooth families of invariant tori}

\author{Mikhail B. Sevryuk}

\ead{2421584@mail.ru, sevryuk@mccme.ru}

\address{V.L.~Tal'rose Institute for Energy Problems of Chemical Physics, N.N.~Sem\"enov Federal Research Center of Chemical Physics, Russian Academy of Sciences, 38 Leninski\u{\i} Prospect, Bld.~2, Moscow 119334, Russia}

\begin{abstract}
For various values of $n$, $d$, and the phase space dimension, we construct simple examples of Hamiltonian and reversible systems possessing smooth $d$-parameter families of invariant $n$-tori carrying conditionally periodic motions. In the Hamiltonian case, these tori can be isotropic, coisotropic, or atropic (neither isotropic nor coisotropic). The cases of non-compact and compact phase spaces are considered. In particular, for any $N\geq 3$ and any vector $\omega\in\mR^N$, we present an example of an analytic Hamiltonian system with $N$ degrees of freedom and with an isolated (and even unique) invariant $N$-torus carrying conditionally periodic motions with frequency vector $\omega$ (but this torus is atropic rather than Lagrangian and the symplectic form is not exact). Examples of isolated atropic invariant tori carrying conditionally periodic motions are given in the paper for the first time. The paper can also be used as an introduction to the problem of the isolatedness of invariant tori in Hamiltonian and reversible systems.
\end{abstract}

\begin{keyword}
Hamiltonian systems; Reversible systems; Kronecker torus; Isolatedness; Uniqueness; Families of tori; Lagrangian torus; Isotropic torus; Coisotropic torus; Atropic torus; Symmetric torus; KAM theory

\MSC[2020] 70K43 \sep 70H12 \sep 70H33 \sep 70H08
\end{keyword}

\end{frontmatter}

\section{Introduction and overview}\label{introduction}

\subsection{Kronecker tori}\label{Kronecker}

Finite-dimensional invariant tori carrying conditionally periodic motions are among the key elements of the structure of smooth dynamical systems with continuous time. The importance and ubiquity of such tori stems, in the long run, from the fact that any finite-dimensional connected compact Abelian Lie group is a torus \cite{A1969,DK2000,S2007}. By definition, given a certain flow on a certain manifold, an invariant $n$-torus carrying \emph{conditionally periodic} motions ($n$ being a non-negative integer) is an invariant submanifold $\cT$ diffeomorphic to the standard $n$-torus $\mT^n=(\mR/2\pi\mZ)^n$ and such that the induced dynamics on $\cT$ in a suitable angular coordinate $\varp\in\mT^n$ has the form $\dot{\varp}=\omega$ where $\omega\in\mR^n$ is a constant vector (called the \emph{frequency vector}). Flows on $\mT^n$ afforded by equations $\dot{\varp}\equiv\omega$ are also said to be linear, parallel, rotational, translational, or Kronecker, and invariant tori carrying conditionally periodic motions are therefore sometimes called \emph{Kronecker tori} \cite{KP2003,MP643,P707,S415,T2012}.

A Kronecker flow $g^t$ on $\mT^n$ with any frequency vector $\omega\in\mR^n$ possesses the \emph{uniform recurrence property}: for any $T>0$ and $\vare>0$ there exists $\Theta\geq T$ such that for any $\varp\in\mT^n$ the distance between $\varp$ and $g^\Theta(\varp) = \varp+\Theta\omega$ (e.g., with respect to some fixed Riemannian metric) is smaller than $\vare$. Recall the almost obvious proof of this fact. First of all, there is $\delta>0$ such that the distance between $\varp$ and $\varp+\Delta$ is smaller than $\vare$ whenever $\varp\in\mT^n$ and $\Delta\in\mR^n$, $|\Delta|<\delta$ (here and henceforth, given $c\in\mR^n$, the symbols $|c|$ denote the $\Fl_1$-norm $|c_1|+\cdots+|c_n|$ of $c$). Second, there are positive integers $\Fm_2>\Fm_1$ and a vector $\Delta\in\mR^n$, $|\Delta|<\delta$ such that $\Fm_2T\omega = \Fm_1T\omega+\Delta$. Finally, it is sufficient to set $\Theta=(\Fm_2-\Fm_1)T$, because $\varp+(\Fm_2-\Fm_1)T\omega = \varp+\Delta$ for any $\varp\in\mT^n$.

The frequency vector $\omega=(\omega_1,\ldots,\omega_n) \in \mR^n$ of a Kronecker $n$-torus and the torus itself are said to be \emph{non-resonant} if the frequencies $\omega_1,\ldots,\omega_n$ are incommensurable (linearly independent over rationals) and are said to be \emph{resonant} otherwise. Conditionally periodic motions with non-resonant frequency vectors are usually called \emph{quasi-periodic} motions. Each phase curve on a non-resonant Kronecker torus fills up it densely. If the frequencies $\omega_1,\ldots,\omega_n$ of a resonant Kronecker $n$-torus $\cT$ satisfy $r$ independent resonance relations $\bigl\langle j^{(\iota)},\omega \bigr\rangle=0$, $j^{(\iota)}\in\mZ^n\setminus\{0\}$, $1\leq\iota\leq r$, $1\leq r\leq n$ (here and henceforth, the angle brackets $\langle{\cdot},{\cdot}\rangle$ denote the standard inner product), then $\cT$ is foliated by non-resonant Kronecker $(n-r)$-tori, the frequency vector of all these tori being the same.

The occurrence of a resonant Kronecker torus in a certain dynamical system usually indicates some degeneracy (for instance, the presence of many first integrals). Normally, dynamical systems exhibit (smooth or Cantor-like) families of non-resonant Kronecker tori, the dimension of all the tori in a given family being the same. Moreover, the frequency vectors of these tori are not merely non-resonant but \emph{strongly non-resonant} (for instance, Diophantine), i.e., badly approximable by resonant vectors. Recall that a vector $\omega\in\mR^n$ is said to be \emph{Diophantine} if there exist constants $\tau\geq n-1$ and $\gamma>0$ such that $\bigl| \langle j,\omega\rangle \bigr| \geq \gamma|j|^{-\tau}$ for any $j\in\mZ^n\setminus\{0\}$ (vectors that are not Diophantine are said to be \emph{Liouville}). Families of Kronecker tori with strongly incommensurable frequencies are the subject of the \emph{KAM} (Kolmogorov--Arnold--Moser) theory. The reader is referred to e.g.\ the monographs \cite{BHS1996,BHTB1990,KP2003}, \S\S~6.2--6.4 of the monograph \cite{AKN2006}, and the survey or tutorial papers \cite{BS249,dlL175,S1113,S137,S603} for the main ideas, methods, and results of the (mainly finite-dimensional) KAM theory and the bibliography, as well as for a precise definition of a Cantor-like family of Kronecker tori. The book \cite{D2014} presents a brilliant semi-popular introduction to the KAM theory. The ``core'' of the theory, namely, families of Kronecker $N$-tori in Hamiltonian systems with $N$ degrees of freedom, is treated in detail in e.g.\ the articles \cite{B42,B21,EFK1733,MP643,P707,S351,T12851}. The papers dealing with various special aspects of the KAM theory are exemplified by \cite{BHN355,BH191,F1521,H79,H989,H49,K259,P380,QS757,S435,S599,S415}. Some open problems in the theory are listed and discussed in the works \cite{FK1905,H797,St177,S6215}.

Typical finite-dimensional autonomous dissipative systems (with no special structure on the phase space the system is assumed to preserve) possess equilibria (Kronecker $0$-tori) and closed trajectories (Kronecker $1$-tori with a nonzero frequency). Typical smooth families of dissipative systems depending on $\Fr\geq 1$ external parameters $\mu_1,\ldots,\mu_{\Fr}$ also exhibit Cantor-like $\Fr$-parameter families of strongly non-resonant Kronecker $n$-tori in the product of the phase space and the parameter space $\{\mu\}$, the dimension $n$ ranging between $2$ and the phase space dimension \cite{BHN355,BHS1996,BHTB1990}. Here and henceforth, the word \emph{``typical''} means that the systems (or the families of systems) with the properties indicated constitute an open set (to be more precise, a set with non-empty interior) in the appropriate functional space.

On the other hand, finite-dimensional autonomous Hamiltonian and reversible systems typically admit many Cantor-like families of strongly non-resonant Kronecker $n$-tori, and these families (with different dimensions $n$) constitute complicated hierarchical structures.

\subsection{Review and the main result:\ Hamiltonian systems}\label{gamiltonovy}

In this section and henceforth, we will employ the following useful notation. Given a non-negative integer $a$, the combination of symbols $\mR^a_w$ will denote the Euclidean space $\mR^a$ with coordinates $(w_1,\ldots,w_a)$, and the combination of symbols $\mT^a_\varp$ will denote the torus $\mT^a$ with angular coordinates $(\varp_1,\ldots,\varp_a)$.

The properties of Kronecker tori in Hamiltonian systems very much depend on the ``relations'' of the torus in question with the symplectic $2$-form. Recall that a submanifold $\cL$ of a $2N$-dimensional symplectic manifold is said to be \emph{isotropic} if the tangent space $T_\Lambda\cL$ to $\cL$ at any point $\Lambda\in\cL$ is contained in its skew-orthogonal complement: $T_\Lambda\cL\subset(T_\Lambda\cL)^\bot$ (in other words, if the restriction of the symplectic form to $\cL$ vanishes), and is said to be \emph{coisotropic} if the tangent space $T_\Lambda\cL$ to $\cL$ at any point $\Lambda\in\cL$ contains its skew-orthogonal complement: $(T_\Lambda\cL)^\bot\subset T_\Lambda\cL$. If $\cL$ is isotropic then $\dim\cL\leq N$, and if $\cL$ is coisotropic then $\dim\cL\geq N$. A submanifold $\cL$ that is both isotropic and coisotropic is said to be \emph{Lagrangian}, in which case $\dim\cL=N$. In the sequel, it will be convenient to call isotropic submanifolds $\cL$ with $\dim\cL<N$ \emph{strictly isotropic} and to call coisotropic submanifolds $\cL$ with $\dim\cL>N$ \emph{strictly coisotropic}. In other words, strictly isotropic submanifolds are isotropic submanifolds that are not Lagrangian, and strictly coisotropic submanifolds are coisotropic submanifolds that are not Lagrangian.

\begin{rem}\label{strictly}
In the literature, the terms ``strictly isotropic'' and ``strictly coisotropic'' are also used with a different meaning, see e.g.\ \cite{B2005,BZ365}. In the $h$-principle theory, one speaks of subcritical isotropic immersions and embeddings in symplectic and contact manifolds where the meaning of the words ``subcritical isotropic'' is close to ``non-Lagrangian isotropic'' (see e.g.\ the tutorial \cite{EM2002}).
\end{rem}

It is clear that if $\dim\cL$ is equal to $0$ or $1$ then $\cL$ is necessarily isotropic, and if $\codim\cL$ is equal to $0$ or $1$ then $\cL$ is necessarily coisotropic. Now suppose that $\cL$ is invariant under a Hamiltonian flow with Hamilton function $H$, $H|_{\cL}$ is a constant, and almost all the points of $\cL$ are not equilibria. Then $\cL$ is isotropic if $\dim\cL=2$ and is coisotropic if $\codim\cL=2$. Recall the simple proof of this fact. Denote the symplectic form by $\Omega$. Let a point $\Lambda\in\cL$ be not an equilibrium, and let $X\in T_\Lambda\cL\setminus\{0\}$ be the vector of the Hamiltonian vector field in question at $\Lambda$. Let $\dim\cL=2$. For any vector $Y\in T_\Lambda\cL$ one has $\Omega(Y,X) = dH(Y) = 0$ since $H$ is a constant on $\cL$. Consequently, $\cL$ is isotropic. On the other hand, let $\codim\cL=2$, and let $\fH\supset\cL$ be the level hypersurface of $H$ containing $\cL$. The space of the tangent vectors to the ambient symplectic manifold at $\Lambda$ that are skew-orthogonal to $X$ is just $T_\Lambda\fH \supset T_\Lambda\cL$, in particular, $X\in(T_\Lambda\cL)^\bot$. Let $X,Y$ be a basis of $(T_\Lambda\cL)^\bot$. Then $Y\in T_\Lambda\fH$ since $X\in T_\Lambda\cL$ and therefore $\Omega(Y,X)=0$. If $Y\in T_\Lambda\fH\setminus T_\Lambda\cL$, then $Y$ would be skew-orthogonal to the whole space $T_\Lambda\fH$ because $\Omega(Y,Y)=0$ and $\dim T_\Lambda\fH-\dim T_\Lambda\cL = 1$, so that $\dim(T_\Lambda\fH)^\bot \geq 2$ in this hypothetical case. Thus, $Y\in T_\Lambda\cL$. Consequently, $\cL$ is coisotropic.

One of the key facts in the Hamiltonian KAM theory is the \emph{Herman lemma} which states that any non-resonant Kronecker torus of a Hamiltonian system is isotropic provided that the symplectic form is exact (see e.g.\ \cite{BHS1996,F1521,S1113} for a proof and \cite{BS249,S137} for a discussion; these works also contain references to the original papers by M.R.~Herman). A Hamiltonian system on a symplectic manifold with a non-exact symplectic form may admit strictly coisotropic non-resonant Kronecker tori as well as non-resonant Kronecker tori that are neither isotropic nor coisotropic. Tori of the latter type are said to be \emph{atropic} \cite{BS249,S1113,St177}. According to what was explained in the previous paragraph, the dimension of an atropic non-resonant Kronecker torus always lies between $3$ and $2N-3$ where $N$ is the number of degrees of freedom.

Now the main ``informal'' conclusion of the Hamiltonian KAM theory can be stated as follows. Typical Hamiltonian systems with $N\geq 1$ degrees of freedom admit $n$-parameter families of isotropic strongly non-resonant Kronecker $n$-tori for each $0\leq n\leq N$. These families are smooth for $n=0$ and $1$ and are Cantor-like for $n\geq 2$. If $N\geq 2$ and the symplectic form is not exact (and meets certain Diophantine-like conditions), typical Hamiltonian systems also possess $(2N-n)$-parameter families of strictly coisotropic strongly non-resonant Kronecker $n$-tori for each $N+1\leq n\leq 2N-1$ (see the works \cite{BHS1996,BS249,H989,H49,P380,S1113,St177} and references therein). These families are smooth for $n=2N-1$ and are Cantor-like for $n\leq 2N-2$. Finally, if $N\geq 3$ and the symplectic form is not exact (and meets certain Diophantine-like conditions), typical Hamiltonian systems also exhibit Cantor-like $\kappa$-parameter families of atropic strongly non-resonant Kronecker $n$-tori for any $3\leq n\leq 2N-3$ and $1\leq\kappa\leq\min(n-2, \, 2N-n-2)$ such that $n+\kappa$ is even (see the works \cite{BS249,S1113,St177} and references therein).

\begin{rem}\label{dichotomyHam}
One sees that if $\kappa$ is the number of parameters in typical families of Kronecker $n$-tori in Hamiltonian systems with $N$ degrees of freedom, then $\kappa+n=2N$ for coisotropic (Lagrangian or strictly coisotropic) Kronecker tori (so that the Lebesgue measure of the union of the tori is positive) and $\kappa+n\leq 2N-2$ for non-coisotropic (strictly isotropic or atropic) Kronecker tori (so that the union of the tori is of measure zero).
\end{rem}

Until 1984, the Hamiltonian KAM theory only dealt with isotropic Kronecker tori. The non-isotropic Hamiltonian KAM theory was founded by I.O.~Parasyuk \cite{P380}. Strictly isotropic Kronecker tori in Hamiltonian systems are often said to be \emph{lower dimensional}. For strictly coisotropic Kronecker tori in Hamiltonian systems, the term ``higher dimensional'' is also used but much more rarely.

Of course, of all the Kronecker $n$-tori in Hamiltonian systems with $N$ degrees of freedom, Lagrangian Kronecker $N$-tori are best studied. A generic Lagrangian non-resonant Kronecker torus $\cT$ in a Hamiltonian system (with any number $N$ of degrees of freedom) is \emph{KAM stable}: in any neighborhood of $\cT$, there is a family of other Lagrangian Kronecker tori, their union having positive Lebesgue measure and density one at $\cT$ (all these tori constitute an $N$-parameter family which is, generally speaking, Cantor-like for $N\geq 2$). To be more precise, the KAM stability of $\cT$ is implied by the so-called Kolmogorov non-degeneracy of $\cT$ \cite{B42}. No arithmetic conditions (like strong incommensurability) on the frequencies of $\cT$ are needed in this remarkable result, and it is valid in the $C^\ell$ smoothness class with any finite sufficiently large $\ell$ (not to mention the $C^\infty$, Gevrey, and analytic categories). On the other hand, generic Lagrangian resonant Kronecker tori in Hamiltonian systems with $N\geq 2$ degrees of freedom in the $C^\ell$ smoothness classes, $\ell\geq 2$, are not KAM stable \cite{B21}. For some previous results concerning density points of quasi-periodicity (not necessarily in the Hamiltonian realm), see e.g.\ the works \cite{BHN355,BHS1996,BHTB1990,BS249,EFK1733}.

\begin{rem}\label{stability}
The term ``KAM stable'' is sometimes understood in a quite different sense (see e.g.\ \cite{K259,S351}): an unperturbed system (or an unperturbed Hamilton function) possessing a smooth family of Kronecker tori is said to be KAM stable if any perturbed system admits a Cantor-like family of Kronecker tori close to the unperturbed ones (provided that the perturbation lies in the suitable functional class and is sufficiently small).
\end{rem}

Since Kronecker tori in Hamiltonian systems tend to be organized into (Cantor-like) families, the natural question arises whether such tori can be isolated. The isolatedness of a torus can be understood in different ways.

\begin{dfn}\label{isolated}
A Kronecker $n$-torus $\cT$ of a dynamical system is said to be \emph{isolated} if it is not included in a (Cantor-like) family of Kronecker $n$-tori. A torus $\cT$ is said to be \emph{strongly isolated} if there exists a neighborhood $\cO$ of $\cT$ in the phase space such that there are no Kronecker tori (of any dimension) entirely contained in $\cO\setminus\cT$. A torus $\cT$ is said to be \emph{unique} if there are no Kronecker tori (of any dimension) outside $\cT$ in the whole phase space.
\end{dfn}

In particular, a strongly isolated torus $\cT$ is unique in the neighborhood $\cO$ mentioned in Definition~\ref{isolated}.

Of course, a generic equilibrium in a Hamiltonian system is always isolated, to be more precise, an equilibrium $O$ is isolated whenever none of the eigenvalues $\lambda_{\Fi}$ of the linearization of the vector field at $O$ is zero. If $O$ is \emph{hyperbolic} (i.e., if all the eigenvalues $\lambda_{\Fi}$ have nonzero real parts), then it is strongly isolated. Examples of unique equilibria in Euclidean phase spaces are also straightforward: if the equilibrium $0$ of a system with a quadratic Hamilton function in $\mR^{2N}$ is hyperbolic then it is unique. On the other hand, the question of whether an \emph{elliptic} equilibrium of a Hamiltonian system (i.e., an equilibrium for which all the eigenvalues $\lambda_{\Fi}$ are nonzero and lie on the imaginary axis) can be strongly isolated (or at least can be not accumulated by a set of Lagrangian Kronecker tori of positive measure) is very far from being easy. So is the question of whether such an equilibrium can be Lyapunov unstable. We will not discuss this problem here and confine ourselves by citing the papers and preprints \cite{F09059,FK1905,FS67,T12851} (see also the references therein). In general, the instability of an equilibrium is a more delicate topic than that of, say, a Lagrangian Kronecker torus \cite{FS67}.

Surprisingly, it seems that the question of whether strictly isotropic Kronecker tori of dimensions from $1$ to $N-1$ in Hamiltonian systems with $N\geq 2$ degrees of freedom can be isolated was never considered until 2017. In December 2017 and January 2018, the author and the user Khanickus of MathOverflow \cite{K2018} constructed independently two very similar explicit (and exceedingly simple) examples of Hamiltonian systems in $\mR^3\times\mT^1$ with a unique periodic orbit. Subsequently, for any integers $n\geq 1$ and $N\geq n+1$ and for any vector $\omega\in\mR^n$, the author \cite{S415} proposed an example of a Hamiltonian system with $N$ degrees of freedom, with the phase space $\mR^{2N-n}_w\times\mT^n_\varp$, with the exact symplectic form
\[
\sum_{i=1}^n dw_i\wedge d\varp_i + \sum_{\nu=1}^{N-n} dw_{n+\nu}\wedge dw_{N+\nu},
\]
and with a Hamilton function independent of $\varp$, polynomial in $w$, and such that $\{w=0\}$ is a unique Kronecker $n$-torus (in the sense of Definition~\ref{isolated}), the frequency vector of $\{w=0\}$ being $\omega$. The paper \cite{S415} also contains an example of a Hamiltonian system with $N$ degrees of freedom, with the compact phase space $\mT^{2N-n}_w\times\mT^n_\varp$, with the symplectic form given by the same formula (but no longer exact), and with a trigonometric polynomial Hamilton function independent of $\varp$ and such that $\{w=0\}$ is a strongly isolated Kronecker $n$-torus (in the sense of Definition~\ref{isolated}), the frequency vector of $\{w=0\}$ being $\omega$. Thus, the problem of the possible isolatedness of strictly isotropic Kronecker tori in Hamiltonian systems has been completely solved by now.

It is clear that periodic orbits (Kronecker $1$-tori with a nonzero frequency) of Hamiltonian systems with one degree of freedom are always included in smooth one-parameter families (each periodic orbit being a connected component of an energy level line). The question of whether Lagrangian Kronecker tori in Hamiltonian systems with $N\geq 2$ degrees of freedom can be isolated has turned out to be highly nontrivial. To the best of the author's knowledge, this question is still open, even in the case where the frequency vector is Liouville (or even resonant) and the Hamilton function is only $C^\infty$ smooth (see \cite{B42,B21,FF01575}).

It is proven in the landmark paper \cite{EFK1733} that a Lagrangian Kronecker $N$-torus $\cT$ with a \emph{Diophantine} frequency vector is never isolated in the \emph{analytic} category (where the symplectic form, the Hamilton function, and the torus itself are analytic), however degenerate the Hamilton function is at $\cT$. Such a torus is always accumulated by other Lagrangian Kronecker tori (with Diophantine frequency vectors), i.e., is always included in a (Cantor-like) $\Fr$-parameter family of Lagrangian Kronecker tori with $\Fr\geq 1$. Nevertheless, it is not known whether $\Fr$ is always equal to $N$, i.e., whether the union of Lagrangian Kronecker tori in any neighborhood of $\cT$ always has positive measure (Herman conjectured the affirmative answer in the model problem of fixed points of analytic symplectomorphisms \cite{H797}). For $N=2$, however, the equality $\Fr=2$ is always valid even in the $C^\infty$ category \cite{EFK1733}.

For Liouville frequency vectors or non-analytic Hamilton functions, there are several examples in the literature of a Lagrangian Kronecker $N$-torus $\cT$ that is accumulated by other Lagrangian Kronecker tori, but the union of these tori is of measure zero. The phase space in all these examples is $\mR^N_u\times\mT^N_\varp$, the symplectic form is $\sum_{i=1}^N du_i\wedge d\varp_i$, the Hamilton function is $H(u,\varp) = \langle u,\omega\rangle + O\bigl( |u|^2 \bigr)$, and $\cT=\{u=0\}$, where $\omega\in\mR^N$ is the frequency vector of $\cT$. It is well known that locally, in some neighborhood of a Lagrangian Kronecker torus, this setup can always be achieved (up to an additive constant in $H$) \cite{B42,B21,MP643}; in fact, this is an immediate consequence of A.~Weinstein's equivalence theorem for Lagrangian submanifolds \cite{W329}. For any $N\geq 2$ and any resonant vector $\omega$, a very simple example with an analytic (and even quadratic in $u$) Hamilton function $H$ is presented in the paper \cite{B21}. In this example, the torus $\cT$ is accumulated by a continuous $N$-parameter family of isotropic Kronecker $(N-1)$-tori. The article \cite{EFK1733} contains examples for any $N\geq 4$ and \emph{any} vector $\omega$ with $C^\infty$ as well as Gevrey regular (with any exponent $\sigma>1$) Hamilton functions $H$. It is pointed out in the paper \cite{FS67} that for $C^\infty$ Hamilton functions $H$, the construction of \cite{EFK1733} can be extended to $N=3$ and any vector $\omega$ (and an analog for elliptic equilibria in $\mR^6$ is described). The paper \cite{FS67} also presents an analog for elliptic equilibria in $\mR^4$ but for Liouville frequencies. Finally, for any $N\geq 3$ and non-resonant but ``sufficiently Liouville'' vectors $\omega$, G.~Farr\'e and B.~Fayad \cite{FF01575} constructed examples with \emph{analytic} Hamilton functions $H$. The words ``sufficiently Liouville'' mean that if $\tilde{\omega} = (\omega_1,\ldots,\omega_{N-1}) \in \mR^{N-1}$, then the infimum of the set of the ratios
\[
\frac{\ln\bigl| \langle j,\tilde{\omega}\rangle \bigr|}{|j|}, \quad j\in\mZ^{N-1}\setminus\{0\}
\]
is $-\infty$. In the examples of \cite{EFK1733} and \cite{FF01575}, the hypersurface $\{u_N=0\}$ is foliated by Lagrangian Kronecker tori with frequency vector $\omega$.

As far as the author knows, the question of whether strictly coisotropic or atropic Kronecker tori in Hamiltonian systems can be isolated has never been raised. The present paper gives an exhaustive answer to this question in the case of \emph{atropic} tori. Here is our main result.

\begin{thm}\label{mainHam}
For any integers $N$, $n$, $d$ in the ranges $N\geq 2$, $1\leq n\leq N-1$, $0\leq d\leq 2N-n$ and for any vector $\omega\in\mR^n$, there exist an exact symplectic form $\Omega$ on the manifold $\cM = \mR^{2N-n}_w\times\mT^n_\varp$ with constant coefficients and a Hamilton function $H:\cM\to\mR$ independent of $\varp$, polynomial in $w$, and such that the corresponding Hamiltonian system on $\cM$ admits a $d$-parameter analytic family of \emph{strictly isotropic} Kronecker $n$-tori of the form $\{w=\const\}$. There are no Kronecker tori (of any dimension) outside this family. The $n$-torus $\{w=0\}$ belongs to this family, its frequency vector is equal to $\omega$, and if $d=0$ then this torus is unique in the sense of Definition~\ref{isolated}. If $d=2N-n$ then the family in question makes up the whole phase space.

For any integers $N$, $n$, $d$ in the ranges $N\geq 3$, $3\leq n\leq 2N-3$, $0\leq d\leq 2N-n$ and for any vector $\omega\in\mR^n$, there exist a non-exact symplectic form $\Omega$ on the manifold $\cM = \mR^{2N-n}_w\times\mT^n_\varp$ with constant coefficients and a Hamilton function $H:\cM\to\mR$ independent of $\varp$, polynomial in $w$, and such that the corresponding Hamiltonian system on $\cM$ admits a $d$-parameter analytic family of \emph{atropic} Kronecker $n$-tori of the form $\{w=\const\}$. There are no Kronecker tori (of any dimension) outside this family. The $n$-torus $\{w=0\}$ belongs to this family and its frequency vector is equal to $\omega$, and if $d=0$ then this torus is unique. If $d=2N-n$ then the family in question makes up the whole phase space.

Similar statements hold \emph{mutatis mutandis} for the compact manifold $\widehat{\cM} = \mT^{2N-n}_w\times\mT^n_\varp$, the modifications being as follows. First, the symplectic form $\Omega$ is not exact in the case of strictly isotropic Kronecker $n$-tori either. Second, the Hamilton function $H$ is now trigonometric polynomial in $w$. Third, it is no longer valid that there are no Kronecker tori (of any dimension) outside the family under consideration. Fourth, if $d=0$ then the $n$-torus $\{w=0\}$ is strongly isolated (rather than unique) in the sense of Definition~\ref{isolated}.
\end{thm}

This theorem is proven in Sections~\ref{symplectic}--\ref{compact} by constructing explicit examples which generalize the examples of the note \cite{S415}. The problem of whether strictly coisotropic Kronecker tori in Hamiltonian systems can be isolated remains open. There is little doubt that this problem is as difficult as the analogous problem (discussed above) for Lagrangian Kronecker tori. Thus, the isolatedness question for Kronecker tori in Hamiltonian systems is very hard for coisotropic (Lagrangian or strictly coisotropic) Kronecker $n$-tori (for $n\geq 2$) and is rather easy for non-coisotropic (strictly isotropic or atropic) ones. This dichotomy surprisingly coincides with the other dichotomy pointed out in Remark~\ref{dichotomyHam}.

Setting $n=N\geq 3$ and $d=0$ in Theorem~\ref{mainHam}, we obtain a Hamilton function $H: \mR^n_w\times\mT^n_\varp \to \mR$ independent of $\varp$, polynomial in $w$, and such that the $n$-torus $\{w=0\}$ is a unique Kronecker torus, the frequency vector of this torus can be any prescribed vector in $\mR^n$. However, this astonishing picture is marred by the fact that the corresponding symplectic form on $\mR^n_w\times\mT^n_\varp$ is not standard and even not exact, and the torus $\{w=0\}$ is atropic rather than Lagrangian.

In the examples of Sections~\ref{symplectic}--\ref{compact}, one deals with coisotropic and non-coisotropic Kronecker tori in a unified way. However, coisotropic Kronecker $n$-tori in our examples for $N\geq 1$ degrees of freedom ($N\leq n\leq 2N-1$) are always organized into $(2N-n)$-parameter analytic families.

Most probably, non-resonant coisotropic Kronecker $(2N-1)$-tori in Hamiltonian systems with $N$ degrees of freedom cannot be isolated for any $N\geq 2$, cf.\ \cite{BHS1996,H989,H49}.

\subsection{Review and the main result:\ reversible systems}\label{obratimye}

While speaking of Kronecker tori (and, more generally, any invariant submanifolds) in reversible systems, one usually only considers \emph{symmetric} invariant submanifolds, i.e., invariant submanifolds that are also invariant under the reversing involution $G$ of the phase space. The dynamics of $G$-reversible systems and the properties of symmetric invariant submanifolds in such systems very much depend on the structure of the fixed point set $\Fix G$ of the involution $G$. This set is a submanifold of the phase space of the same smoothness class as the involution $G$ itself. However, the manifold $\Fix G$ can well be empty or consist of several connected components of different dimensions even if the phase space is connected (see e.g.\ simple examples in the papers \cite{BDP223,DP3119,PF280,QS757} and references therein; in fact, the literature on the structure of the fixed point sets of involutions of various manifolds is by now immense).

It is well known that in any symmetric non-resonant Kronecker $n$-torus $\cT$ (with a frequency vector $\omega\in\mR^n$) of a $G$-reversible system, one can choose an angular coordinate $\varp\in\mT^n$ such that the dynamics on $\cT$ takes the form $\dot{\varp}=\omega$ and the restriction of the reversing involution $G$ to $\cT$ takes the form $G|_{\cT}: \varp\mapsto-\varp$ (in particular, this implies that the set $(\Fix G)\cap\cT = \Fix\bigl( G|_{\cT} \bigr)$ consists of $2^n$ points). This very easy but fundamental \emph{standard reflection lemma} is proven in e.g.\ the works \cite{BHS1996,S435} (see also the papers \cite{S137,S599} for a discussion).

We will confine ourselves with the case where the fixed point set $\Fix G$ of the reversing involution $G$ is non-empty and all its connected components are of the same dimension, so that $\dim\Fix G$ is well defined (in fact, this is so for almost all the reversible systems encountered in practice). We will say that an involution $G$ satisfying this condition is of \emph{type $(\fL,m)$} if $\dim\Fix G=m$ and $\codim\Fix G=\fL$ (cf.\ \cite{BHN355,BH191,BHS1996,QS757}). It follows from the standard reflection lemma that if a system reversible with respect to an involution of type $(\fL,m)$ admits a symmetric non-resonant Kronecker $n$-torus then $n\leq\fL$. Therefore, in the reversible KAM theory \cite{BHN355,BH191,BHS1996,QS757,St177,S435,S137,S599,S603,S415}, it only makes sense to consider symmetric Kronecker $n$-tori in systems reversible with respect to involutions of types $(n+l,m)$ with $l\geq 0$.

Now the main ``informal'' conclusion of the KAM theory for \emph{individual} reversible systems (not for reversible systems depending on external parameters) can be stated as follows. For $n=0$ and $1$, typical systems reversible with respect to an involution of type $(n+l,m)$ with $m\geq l\geq 0$ admit smooth $(m-l)$-parameter families of symmetric non-resonant Kronecker $n$-tori. For each $n\geq 2$, typical systems reversible with respect to an involution of type $(n+l,m)$ with $l\geq 0$ and $m\geq l+1$ admit Cantor-like $(m-l)$-parameter families of symmetric strongly non-resonant Kronecker $n$-tori.

\begin{rem}\label{dichotomyrev}
One sees that if $\kappa=m-l$ is the number of parameters in typical families of symmetric Kronecker $n$-tori in systems reversible with respect to involutions of types $(n+l,m)$, then $\kappa+n=m+n+l$ (i.e., $\kappa=m+l$) for $l=0$ (so that the Lebesgue measure of the union of the tori is positive) and $\kappa+n<m+n+l$ for $l\geq 1$ (so that the union of the tori is of measure zero). Of course, here we suppose that $m\geq l$ for $n\leq 1$ and $m>l$ for $n\geq 2$. If $m$ and $l$ do not meet these conditions, one needs external parameters $\mu_{\Fj}$ to obtain persistent families of symmetric Kronecker $n$-tori. However, the measure of the union of the tori in the product of the phase space and the parameter space $\{\mu\}$ is typically positive for $l=0$ and zero for $l\geq 1$ in this case as well (see \cite{S435,S137,S599,S603}).
\end{rem}

By analogy with Hamiltonian systems, one may ask under what conditions symmetric Kronecker tori in reversible systems can be isolated or unique. When we speak of the isolatedness, strong isolatedness, and unicity of such tori, we still interpret these concepts in strict accordance with Definition~\ref{isolated}: we have in view the absence of other Kronecker tori whatsoever, and not just the absence of other symmetric Kronecker tori.

For any $n\geq 0$ and $l\geq m\geq 0$, it is very easy to construct a system that is reversible with respect to an involution of type $(n+l,m)$ and admits a unique symmetric Kronecker $n$-torus with any prescribed frequency vector $\omega\in\mR^n$. Indeed, let the phase space be $\mR^m_u\times\mT^n_\varp\times\mR^m_v\times\mR^{l-m}_q$ and the reversing involution be
\[
G: (u,\varp,v,q) \mapsto (u,-\varp,-v,-q).
\]
Then $G$ is an involution of type $(n+l,m)$. The system
\[
\dot{u}=v, \quad \dot{\varp}=\omega, \quad \dot{v}=u, \quad \dot{q}_\nu=q_\nu^2
\]
($1\leq\nu\leq l-m$) is reversible with respect to $G$, and $\{u=0, \; v=0, \; q=0\}$ is a unique Kronecker $n$-torus of this system. This Kronecker torus is symmetric, and its frequency vector is $\omega$.

For any integers $n\geq 0$, $m\geq 0$, $l\geq 1$ and for any vector $\omega\in\mR^n$, the note \cite{S415} considers the manifold
\begin{equation}
\cK = \mR^m_u\times\mT^n_\varp\times\mR^l_q
\label{cK}
\end{equation}
equipped with the involution
\begin{equation}
G: (u,\varp,q) \mapsto (u,-\varp,-q)
\label{inv}
\end{equation}
of type $(n+l,m)$ and presents an example of a $G$-reversible system on $\cK$ with the right-hand side independent of $\varp$, polynomial in $(u,q)$, and such that $\{u=0, \; q=0\}$ is a unique symmetric Kronecker $n$-torus, its frequency vector being $\omega$. For the compact manifold
\begin{equation}
\widehat{\cK} = \mT^m_u\times\mT^n_\varp\times\mT^l_q
\label{hatcK}
\end{equation}
equipped with the involution $G$ given by the same formula \eqref{inv} and having the same type, the paper \cite{S415} contains an example of a $G$-reversible system on $\widehat{\cK}$ with the right-hand side independent of $\varp$, trigonometric polynomial in $(u,q)$, and such that $\{u=0, \; q=0\}$ is a strongly isolated symmetric Kronecker $n$-torus, its frequency vector being $\omega$.

In the present paper, we generalize these examples of the note \cite{S415}. Here is our second result.

\begin{thm}\label{mainrev}
For any integers $n$, $m$, $l$, $d_\ast$, $d$ in the ranges $n\geq 0$, $m\geq 0$, $l\geq 1$, $0\leq d_\ast\leq m$, $d_\ast\leq d\leq d_\ast+l$ and for any vector $\omega\in\mR^n$, there exists a system of ordinary differential equations on \eqref{cK} reversible with respect to the involution \eqref{inv} of type $(n+l,m)$ and possessing the following properties. The right-hand side of this system is independent of $\varp$ and polynomial in $(u,q)$. The system admits a $d$-parameter analytic family of Kronecker $n$-tori of the form $\{u=\const, \; q=\const\}$. There are no Kronecker tori (of any dimension) outside this family. The family includes a $d_\ast$-parameter analytic subfamily of symmetric Kronecker $n$-tori of the form $\{u=\const, \; q=0\}$. The $n$-torus $\{u=0, \; q=0\}$ belongs to this subfamily, its frequency vector is equal to $\omega$, and if $d_\ast=d=0$ then this torus is unique in the sense of Definition~\ref{isolated}. If $d_\ast=m$ and $d=m+l$ then the $d$-parameter family in question makes up the whole phase space.

Similar statements hold \emph{mutatis mutandis} for the compact phase space \eqref{hatcK}, the modifications being as follows. First, the right-hand side of the system is now trigonometric polynomial in $(u,q)$. Second, it is no longer valid that there are no Kronecker tori (of any dimension) outside the family under consideration. Third, the symmetric Kronecker $n$-tori have the form $\{u=\const, \; q=q^0\}$, where each component of $q^0$ is equal to either $0$ or $\pi$. Fourth, if $d_\ast=d=0$ then the $n$-torus $\{u=0, \; q=0\}$ is strongly isolated (rather than unique) in the sense of Definition~\ref{isolated}.
\end{thm}

This theorem is proven in Sections~\ref{reversible}--\ref{comprev}. In the examples of Sections~\ref{reversible}--\ref{comprev}, one deals with the cases $l=0$ and $l\geq 1$ in a unified way. However, for $l=0$ the whole phase space $\mR^m\times\mT^n$ or $\mT^m\times\mT^n$ in our examples is foliated by symmetric Kronecker $n$-tori.

The only case not covered by the examples above is that of symmetric Kronecker $n$-tori in systems reversible with respect to involutions of types $(n,m)$ with $m\geq 1$. If $n=0$ or $1$ then symmetric Kronecker $n$-tori in such systems are always organized into smooth $m$-parameter families and cannot be isolated \cite{S415}. To the best of the author's knowledge, the question of whether symmetric Kronecker $n$-tori in systems reversible with respect to involutions of types $(n,m)$ can be isolated for $m\geq 1$ and $n\geq 2$ has never been raised and is open. One may conjecture that such tori with Diophantine frequency vectors are never isolated in the analytic category (where the involution, the vector field, and the torus itself are analytic), similarly to Lagrangian Kronecker tori in Hamiltonian systems \cite{EFK1733} (see Section~\ref{gamiltonovy}). Most probably, this question is very hard.

To summarize, the problem of the possible isolatedness of symmetric Kronecker $n$-tori in systems reversible with respect to involutions of types $(n+l,m)$ is rather easy (and has been solved) for $l\geq 1$ and is probably highly nontrivial for $l=0$ (if $n\geq 2$ and $m\geq 1$). Like in the Hamiltonian realm, this dichotomy coincides with the dichotomy of Remark~\ref{dichotomyrev}.

\section{Preliminaries}\label{symplectic}

Given non-negative integers $a$ and $b$, we designate the identity $a\times a$ matrix as $I_a$ and the zero $a\times b$ matrix as $0_{a\times b}$. In fact, the symbols $I_0$, $0_{0\times b}$, and $0_{a\times 0}$ correspond to no actual objects and will only be used for unifying the notation.

Let $s\geq 1$, $k$, and $l$ be non-negative integers and consider a skew-symmetric $(2s+2k)\times (2s+2k)$ matrix $J$ of the form
\[
J = \begin{pmatrix}
0_{s\times s} & -Z^{\rt} \\ Z & L
\end{pmatrix},
\]
where $Z$ is an $(s+2k)\times s$ matrix of rank $s$ and $L$ is a skew-symmetric $(s+2k)\times (s+2k)$ matrix (the superscript ``t'' denotes transposing). If $k=0$ then the matrix $J$ is always non-singular ($\det J = (\det Z)^2$). If $k\geq 1$ then for any fixed matrix $Z$, the matrix $J$ may be non-singular or singular depending on the matrix $L$. Indeed, we can suppose without loss of generality that the last $s$ rows of $Z$ constitute a non-singular $s\times s$ matrix $Z_\sharp$. Let the matrix $L$ have the form
\[
L = \begin{pmatrix}
L_\sharp & 0_{2k\times s} \\ 0_{s\times 2k} & 0_{s\times s}
\end{pmatrix},
\]
where $L_\sharp$ is a skew-symmetric $2k\times 2k$ matrix, then $\det J = \det L_\sharp(\det Z_\sharp)^2 \neq 0$ if and only if $\det L_\sharp \neq 0$.

In the sequel, we will assume the matrix $J$ to be non-singular, so that the skew-symmetric $(2s+2k+2l)\times (2s+2k+2l)$ matrix
\begin{equation}
\cJ = \begin{pmatrix}
\begin{matrix} 0_{s\times s} & -Z^{\rt} \\ Z & L \end{matrix}
& 0_{(2s+2k)\times 2l} \\ 0_{2l\times (2s+2k)} &
\begin{matrix} 0_{l\times l} & -I_l \\ I_l & 0_{l\times l} \end{matrix}
\end{pmatrix}
\label{cJ}
\end{equation}
is also non-singular and can be treated as the \emph{structure matrix} (the matrix of the Poisson brackets $\{{\cdot},{\cdot}\}$ of the coordinate functions, see e.g.\ \cite{B1,HLW2006,O1993,St177}) of a certain symplectic form $\Omega$ (with constant coefficients) on the manifold
\begin{equation}
\cM = \mR^s_u\times\mT^{s+2k}_\varp\times\mR^l_p\times\mR^l_q
\label{cM}
\end{equation}
(cf.\ Lemma~1 in \cite{St177}). A Hamilton function $H:\cM\to\mR$ affords the equations of motion \cite{B1,HLW2006,O1993}
\begin{equation}
\begin{pmatrix}
\dot{u} \\ \dot{\varp} \\ \dot{p} \\ \dot{q}
\end{pmatrix}
= \cJ\frac{\partial H}{\partial(u,\varp,p,q)} =
\begin{pmatrix}
-Z^{\rt}\partial H/\partial\varp \\
Z\partial H/\partial u + L\partial H/\partial\varp \\
-\partial H/\partial q \\
\partial H/\partial p
\end{pmatrix}.
\label{XH}
\end{equation}
This is an autonomous Hamiltonian system with $N=s+k+l$ degrees of freedom.

For any $u^0\in\mR^s$, $p^0\in\mR^l$, $q^0\in\mR^l$, consider the $(s+2k)$-torus
\begin{equation}
\cT_{u^0,p^0,q^0} = \bigl\{ (u^0,\varp,p^0,q^0) \bigm| \varp\in\mT^{s+2k} \bigr\}.
\label{cT}
\end{equation}
For any $\varp^0\in\mT^{s+2k}$, the skew-orthogonal complement $T^\bot$ (with respect to $\Omega$) of the tangent space $T$ to $\cT_{u^0,p^0,q^0}$ at the point $(u^0,\varp^0,p^0,q^0)$ consists of all the vectors of the form $\psi\partial/\partial\varp + P\partial/\partial p + Q\partial/\partial q$, where $P\in\mR^l$, $Q\in\mR^l$, $\psi\in\cZ$, and $\cZ$ is the $s$-dimensional subspace of $\mR^{s+2k}$ spanned by the columns of the matrix $Z$. Indeed, the space of all such vectors is of dimension $s+2l = \dim\cM-(s+2k)$. It is therefore sufficient to verify that $\Omega(V,W) = 0$ for any vector $V\in T$ (i.e., any vector $V = \Phi\partial/\partial\varp$ with $\Phi\in\mR^{s+2k}$) and any vector $W = (ZU)\partial/\partial\varp + P\partial/\partial p + Q\partial/\partial q$ with $U\in\mR^s$, $P\in\mR^l$, $Q\in\mR^l$. According to \eqref{XH}, the linear Hamilton function $H = H(u,p,q) = \langle U,u\rangle-\langle P,q\rangle+\langle Q,p\rangle$ on $\cM$ affords the constant Hamiltonian vector field equal to $W$. Thus, $\Omega(V,W) = dH(V) = 0$.

We arrive at the conclusion that the $(s+2k)$-tori \eqref{cT} are isotropic for $k=0$ ($T\subset T^\bot$ at any point), are coisotropic for $l=0$ ($T^\bot\subset T$ at any point), and are therefore Lagrangian for $k=l=0$. For $kl>0$, these tori are atropic. Note that $\dim(T\cap T^\bot)=s$ in all the cases.

It is clear that the symplectic form $\Omega$ is exact if and only if its coordinate representation does not contain terms $c_{\alpha\beta}d\varp_\alpha\wedge d\varp_\beta$, $1\leq\alpha<\beta\leq s+2k$, i.e., if the tori \eqref{cT} are isotropic. Thus, $\Omega$ is exact for $k=0$ and is not exact for $k\geq 1$.

\section{The main construction}\label{center}

\subsection{The system}\label{system}

Now let $\zeta_1,\ldots,\zeta_s$, $\xi_1,\ldots,\xi_l$, $\eta_1,\ldots,\eta_l$ be arbitrary non-negative real constants and let $h:\mR^s\to\mR$ be an arbitrary smooth function. Consider the Hamilton function
\begin{equation}
H(u,p,q) = h(u) + lp_1 \sum_{i=1}^s \zeta_iu_i^2 + \sum_{\nu=1}^l (\xi_\nu p_\nu q_\nu^2+\eta_\nu p_\nu^3/3)
\label{ourH}
\end{equation}
on the symplectic manifold \eqref{cM}. The term $lp_1 \sum_{i=1}^s \zeta_iu_i^2$ is automatically absent for $l=0$. According to \eqref{XH}, the equations of motion afforded by $H$ take the form
\begin{equation}
\begin{aligned}
\dot{u}_i &= 0, \\
\dot{\varp}_\alpha &= \sum_{i=1}^s Z_{\alpha i}\left( \frac{\partial h(u)}{\partial u_i} + 2l\zeta_iu_ip_1 \right), \\
\dot{p}_\nu &= -2\xi_\nu p_\nu q_\nu, \\
\dot{q}_\nu &= \xi_\nu q_\nu^2 + \eta_\nu p_\nu^2 + \delta_{1\nu}l \sum_{i=1}^s \zeta_iu_i^2,
\end{aligned}
\label{XourH}
\end{equation}
where $1\leq i\leq s$, $1\leq\alpha\leq s+2k$, $1\leq\nu\leq l$, and $\delta_{1\nu}$ is the Kronecker delta. The fundamental property of this system is that $\dot{q}_\nu\geq 0$ everywhere in the phase space $\cM$, $1\leq\nu\leq l$.

In the note \cite{S415}, we considered the particular case of the Hamilton function \eqref{ourH} and the system \eqref{XourH} where $k=0$, $Z=I_s$, $L=0_{s\times s}$, $h(u)=\langle u,\omega\rangle$ ($\omega\in\mR^s$), $l\geq 1$, $\zeta_i=1/l$ for all $1\leq i\leq s$, and $\xi_\nu=\eta_\nu=1$ for all $1\leq\nu\leq l$ (in the notation of \cite{S415}, $s=n$ and $l=m+1$ where $n\geq 1$ and $m\geq 0$).

All the conditionally periodic motions of the system \eqref{XourH} fill up the manifold
\[
\fM = \bigl\{ (u,\varp,p,q) \bigm| l\zeta_iu_i=0 \; \forall i, \;\; \eta_\nu p_\nu=0 \; \forall\nu, \;\; \xi_\nu q_\nu=0 \; \forall\nu \bigr\}
\]
foliated by Kronecker $(s+2k)$-tori of the form \eqref{cT}. Of course, always $\cT_{0,0,0}\subset\fM$. The frequency vector of a torus $\cT_{u^0,p^0,q^0} \subset \fM$ is $\omega(u^0) = Z\partial h(u^0)/\partial u \in \cZ$ (recall that $\cZ$ is the $s$-dimensional subspace of $\mR^{s+2k}$ spanned by the columns of the matrix $Z$). If $l=0$ then $\fM=\cM$. The system \eqref{XourH} admits no conditionally periodic motions outside $\fM$.

Indeed, if $(u,\varp,p,q)\in\fM$ then $\dot{u}=0$, $\dot{\varp}=Z\partial h(u)/\partial u$, $\dot{p}=0$, $\dot{q}=0$. On the other hand, since $\dot{q}_\nu\geq 0$ everywhere in $\cM$, the recurrence property of conditionally periodic motions implies that $\dot{q}_\nu\equiv 0$ on Kronecker tori, $1\leq\nu\leq l$. Consequently, a point $(u,\varp,p,q)\notin\fM$ does not belong to any Kronecker torus of \eqref{XourH} (of any dimension) because $\dot{q}_1>0$ whenever $l\zeta_iu_i\neq 0$ for at least one $i$ and $\dot{q}_\nu>0$ whenever $\eta_\nu p_\nu\neq 0$ or $\xi_\nu q_\nu\neq 0$, $1\leq\nu\leq l$.

The Kronecker $(s+2k)$-tori $\cT_{u^0,p^0,q^0} \subset \fM$ constitute an analytic $d$-parameter family where $d = \dim\fM-(s+2k)$. If $l=0$ then $d=s$. If $l\geq 1$ (i.e., if the tori \eqref{cT} are not coisotropic) then $d$ can take any integer value between $0$ and $s+2l$; to be more precise, $d$ is the number of zero constants among $\zeta_i$, $\xi_\nu$, $\eta_\nu$ ($1\leq i\leq s$, $1\leq\nu\leq l$). The equality $d=0$ holds if and only if all the numbers $\zeta_i$, $\xi_\nu$, $\eta_\nu$ are positive in which case $\fM=\cT_{0,0,0}$, and $\cT_{0,0,0}$ is a unique Kronecker torus of the system \eqref{XourH}. The equality $d=s+2l$ occurs if and only if all the numbers $\zeta_i$, $\xi_\nu$, $\eta_\nu$ are equal to zero in which case $\fM=\cM$. If a torus $\cT_{u^0,p^0,q^0}$ lies in $\fM$ then its frequency vector $\omega(u^0) = Z\partial h(u^0)/\partial u$ can be made equal to any prescribed vector in $\mR^{s+2k}$ by a suitable choice of the matrix $Z$ and the function $h$ (one can even choose $h$ to be linear).

The construction just described can be formally carried out for $s=0$ as well, but for $s=0$ the frequency vector of each invariant $2k$-torus of the form \eqref{cT} is zero: such a torus consists of equilibria.

The $s+l+\delta_{0l}$ functions
\[
H, \qquad
u_i \;\; (1\leq i\leq s), \qquad
\xi_\nu p_\nu q_\nu^2+\eta_\nu p_\nu^3/3 \;\; (2\leq\nu\leq l)
\]
are first integrals of the system \eqref{XourH} which are pairwise in involution. In fact, this system always admits $s+l$ first integrals that are pairwise in involution and are functionally independent almost everywhere. Indeed, let $f_\nu(p_\nu,q_\nu) = \xi_\nu p_\nu q_\nu^2+\eta_\nu p_\nu^3/3$ if $\xi_\nu+\eta_\nu>0$, and let $f_\nu$ be any smooth function in $p_\nu,q_\nu$ with the differential other than zero almost everywhere if $\xi_\nu=\eta_\nu=0$. In the case where $l\geq 1$ and $\sum_{i=0}^s \zeta_i > 0$, the functions
\[
H, \qquad u_i \;\; (1\leq i\leq s), \qquad f_\nu(p_\nu,q_\nu) \;\; (2\leq\nu\leq l)
\]
are the desired $s+l$ first integrals. In the opposite case where $l \sum_{i=0}^s \zeta_i = 0$, one can choose the $s+l$ first integrals in question to be equal to
\[
u_i \;\; (1\leq i\leq s), \qquad f_\nu(p_\nu,q_\nu) \;\; (1\leq\nu\leq l).
\]

\subsection{The analysis}\label{analysis}

The dimension $n=s+2k$ of the tori \eqref{cT} can be smaller than, equal to, or greater than the number $N=s+k+l$ of degrees of freedom: $n-N=k-l$. The maximal possible value $s+2l$ of the quantity $d$ is always equal to $2N-n$. If $l=0$ (the case of coisotropic tori \eqref{cT}) then $d=s=2N-n$ and $k=n-N$. These equalities determine integers $s\geq 1$ and $k\geq 0$ if and only if $N\geq 1$ and $N\leq n\leq 2N-1$.

If $l\geq 1$ then $n=s+2k=2N-s-2l\leq 2N-3$ and $N\geq 2$. It is easy to see that for any integers $N\geq 2$ and $n$ in the range $1\leq n\leq 2N-3$, one can choose integers $s\geq 1$, $k\geq 0$, $l\geq 1$ such that $N=s+k+l$ and $n=s+2k$. Indeed, if $1\leq n\leq N-1$ then it suffices to set $s=n$, $k=0$, $l=N-n$. In this case, the tori \eqref{cT} are strictly isotropic. Of course, the converse is also true: if $l\geq 1$ and $k=0$ then $1\leq n=s\leq N-1=s+l-1$. On the other hand, if $N\leq n\leq 2N-3$ (so that $N\geq 3$) then it suffices to set $s=2N-n-2$, $k=n-N+1$, $l=1$. In this case, the tori \eqref{cT} are atropic.

If $l\geq 1$ and $k\geq 1$ (so that the tori \eqref{cT} are atropic) then $n=s+2k\geq 3$ and $N=s+k+l\geq 3$. One easily sees that for any integers $N\geq 3$ and $n$ in the range $3\leq n\leq 2N-3$, one can choose positive integers $s$, $k$, $l$ such that $N=s+k+l$ and $n=s+2k$. In the previous paragraph, we verified this for $N\leq n\leq 2N-3$. On the other hand, if $3\leq n\leq N-1$ (so that $N\geq 4$) then it suffices to set $s=n-2$, $k=1$, $l=N-n+1$.

The case where $n=N\geq 3$, $1\leq k=l\leq\bigl\lfloor (N-1)/2 \bigr\rfloor$ (here $\lfloor{\cdot}\rfloor$ denotes the floor function), $s=N-2l$, and $d=0$ is probably the most interesting one. In this case we obtain the unique Kronecker $N$-torus $\cT_{0,0,0}$ of the Hamiltonian system \eqref{XourH} with $N$ degrees of freedom, and the frequency vector of this torus can be any vector in $\mR^N$ (but this torus is atropic rather than Lagrangian).

The case where $n=N\geq 3$, $1\leq k=l\leq\bigl\lfloor (N-1)/2 \bigr\rfloor$, $s=N-2l$, and $d=N$ is also very interesting. In this case the whole phase space of the Hamiltonian system \eqref{XourH} with $N$ degrees of freedom is smoothly foliated by Kronecker $N$-tori \eqref{cT}. However, this system is not \emph{Liouville integrable} (completely integrable): the Kronecker tori in question are atropic rather than Lagrangian, and the $N$ first integrals $u_1,\ldots,u_{N-2l}$, $p_1,\ldots,p_l$, $q_1,\ldots,q_l$ are \emph{not} pairwise in involution: $\{q_\nu,p_\nu\}\equiv 1$, $1\leq\nu\leq l$.

In fact, for any $s\geq 1$, $k$, and $l$, the whole phase space $\cM$ of the Hamiltonian system \eqref{XourH} is smoothly foliated by Kronecker tori \eqref{cT} whenever $d=s+2l$ (as was already pointed out in Section~\ref{system}), in which case the $s+2l$ functions $u_1,\ldots,u_s$, $p_1,\ldots,p_l$, $q_1,\ldots,q_l$ are independent first integrals of the system. The matrix of the Poisson brackets of these functions is
\[
\cP = \begin{pmatrix}
0_{s\times s} & 0_{s\times l} & 0_{s\times l} \\
0_{l\times s} & 0_{l\times l} & -I_l \\
0_{l\times s} & I_l & 0_{l\times l}
\end{pmatrix},
\]
and they are not pairwise in involution for $l\geq 1$. If $l>k$ then the number of the first integrals in question exceeds the number $N=s+k+l$ of degrees of freedom. However, one cannot call the system \eqref{XourH} superintegrable for $l>k\geq 1$ and $d=s+2l$. Besides the existence of $M>N$ independent first integrals, the definition of a \emph{superintegrable} Hamiltonian system with $N$ degrees of freedom (see the papers \cite{F93,H79,KS811} and references therein; superintegrable systems are also known as properly degenerate or non-commutatively integrable systems) includes other requirements, for instance, that there be $N$ integrals pairwise in involution among the $M$ integrals under consideration (we have only $s+l<N$ integrals in involution, e.g., $u_1,\ldots,u_s$, $p_1,\ldots,p_l$), or that the rank of the matrix of the Poisson brackets of the integrals be equal to $2(M-N)$ almost everywhere (in our case the rank of $\cP$ is $2l>2(s+2l-N)=2(l-k)$), or that the common level surfaces of the integrals be isotropic (in our case the tori \eqref{cT} are atropic).

Note that the tori \eqref{cT} are isotropic if and only if the symplectic form $\Omega$ on $\cM$ is exact (both the properties in our setup are equivalent to the equality $k=0$). This observation is consistent with the Herman lemma (see Section~\ref{gamiltonovy}).

\section{Compact phase spaces}\label{compact}

Like in the setting of our note \cite{S415}, the general construction of Sections~\ref{symplectic} and~\ref{center} admits an analogue with a compact phase space. Consider the symplectic manifold
\begin{equation}
\widehat{\cM} = \mT^s_u\times\mT^{s+2k}_\varp\times\mT^l_p\times\mT^l_q
\label{hatcM}
\end{equation}
with the same structure matrix \eqref{cJ}. Of course, now the corresponding symplectic form $\Omega$ is always non-exact. The $(s+2k)$-tori \eqref{cT} with $u^0\in\mT^s$, $p^0\in\mT^l$, $q^0\in\mT^l$ are again isotropic for $k=0$, are coisotropic for $l=0$, and are atropic for $kl>0$. For any angular variable $z$ introduce the notation $\tilde{z}=\sin z$ (cf.\ \cite{S415}). Consider the Hamilton function
\[
\widehat{H}(u,p,q) = h(u) + l\tilde{p}_1 \sum_{i=1}^s \zeta_i\tilde{u}_i^2 + \sum_{\nu=1}^l (\xi_\nu\tilde{p}_\nu\tilde{q}_\nu^2+\eta_\nu\tilde{p}_\nu^3/3)
\]
on \eqref{hatcM}, where again $\zeta_1,\ldots,\zeta_s$, $\xi_1,\ldots,\xi_l$, $\eta_1,\ldots,\eta_l$ are arbitrary non-negative real constants and $h:\mT^s\to\mR$ is an arbitrary smooth function. The Hamilton function $\widehat{H}$ affords the equations of motion
\begin{equation}
\begin{aligned}
\dot{u}_i &= 0, \\
\dot{\varp}_\alpha &= \sum_{i=1}^s Z_{\alpha i}\left( \frac{\partial h(u)}{\partial u_i} + l\zeta_i\sin 2u_i\tilde{p}_1 \right), \\
\dot{p}_\nu &= -\xi_\nu\tilde{p}_\nu\sin 2q_\nu, \\
\dot{q}_\nu &= (\xi_\nu\tilde{q}_\nu^2 + \eta_\nu\tilde{p}_\nu^2)\cos p_\nu + \delta_{1\nu}l\left( \sum_{i=1}^s \zeta_i\tilde{u}_i^2 \right)\cos p_1,
\end{aligned}
\label{XourhatH}
\end{equation}
where $1\leq i\leq s$, $1\leq\alpha\leq s+2k$, $1\leq\nu\leq l$.

The manifold
\[
\widehat{\fM} = \bigl\{ (u,\varp,p,q) \bigm| l\zeta_i\tilde{u}_i=0 \; \forall i, \;\; \eta_\nu\tilde{p}_\nu=0 \; \forall\nu, \;\; \xi_\nu\tilde{q}_\nu=0 \; \forall\nu \bigr\}
\]
is again foliated by Kronecker $(s+2k)$-tori of the form \eqref{cT} (with $u^0\in\mT^s$, $p^0\in\mT^l$, $q^0\in\mT^l$), $\cT_{0,0,0}\subset\widehat{\fM}$ in all the cases, and the frequency vector of a torus $\cT_{u^0,p^0,q^0} \subset \widehat{\fM}$ is $\omega(u^0) = Z\partial h(u^0)/\partial u \in \cZ$. This frequency vector can again be made equal to any prescribed vector in $\mR^{s+2k}$ by a suitable choice of the matrix $Z$ and the function $h$; one can choose $h$ to be of the form $\sum_{i=1}^s c_i\sin(u_i-u^0_i)$. The dimension $s+2k+d$ of the manifold $\widehat{\fM}$ is determined in exactly the same way as that of the manifold $\fM$ in Section~\ref{system}. In particular, if $l=0$ then $\widehat{\fM}=\widehat{\cM}$.

In contrast to the case of the system \eqref{XourH}, it is, generally speaking, \emph{not} true that the system \eqref{XourhatH} for $l\geq 1$ possesses no conditionally periodic motions outside $\widehat{\fM}$. Indeed, suppose that $\sum_{i=1}^s \zeta_i > 0$ and choose an arbitrary point $u^0\in\mT^s$ such that $\chi = \sum_{i=1}^s \zeta_i\sin^2u^0_i > 0$. Consider the $(s+2k+1)$-torus
\begin{equation}
\bigl\{ (u^0,\varp,0,q) \bigm| q_2=\cdots=q_l=0 \bigr\} \not\subset \widehat{\fM}.
\label{exception}
\end{equation}
This torus is invariant under the flow of \eqref{XourhatH} with the induced dynamics
\[
\dot{\varp} = \omega(u^0), \qquad \dot{q}_1 = \xi_1\sin^2q_1 + l\chi.
\]
It is clear that the motion on the torus \eqref{exception} is conditionally periodic. The frequencies of this motion are equal to $\omega_1(u^0),\ldots,\omega_{s+2k}(u^0),\varpi$ where
\[
\varpi = 2\pi\left( \int_0^{2\pi}\frac{d\Fq}{\xi_1\sin^2\Fq+l\chi} \right)^{-1} = \bigl[ l\chi(l\chi+\xi_1) \bigr]^{1/2}.
\]

Nevertheless, for any fixed $q^\star\in\mT^l$, no point $(u,\varp,p,q)\notin\widehat{\fM}$ belongs to a Kronecker torus of \eqref{XourhatH} (of any dimension) entirely contained in the domain
\[
\fD^+_{q^\star} = \bigl\{ (u,\varp,p,q) \bigm| p_\nu\in(-\pi/2,\pi/2)\bmod 2\pi \; \forall\nu, \;\; q_\nu\neq q^\star_\nu \; \forall\nu \bigr\}.
\]
Indeed, $\dot{q}_\nu\geq 0$ everywhere in the domain $\fD^+_{q^\star}$, $1\leq\nu\leq l$, and $\sum_{\nu=1}^l \dot{q}_\nu > 0$ everywhere in $\fD^+_{q^\star} \setminus \widehat{\fM}$. If for some $\nu$ a function $q_\nu: \mR \to \mT^1\setminus\{q^\star_\nu\}$ satisfies the conditions that $\dot{q}_\nu(t)\geq 0$ for all $t$ and $\dot{q}_\nu(0)>0$, then $q_\nu(t)$ tends to a certain point $q_\nu^{\lim}\neq q_\nu(0)$ as $t\to+\infty$ and the recurrence property fails.

Similarly, no point $(u,\varp,p,q)\notin\widehat{\fM}$ belongs to a Kronecker torus of \eqref{XourhatH} (of any dimension) entirely contained in the domain
\[
\fD^-_{q^\star} = \bigl\{ (u,\varp,p,q) \bigm| p_\nu\in(\pi/2,3\pi/2)\bmod 2\pi \; \forall\nu, \;\; q_\nu\neq q^\star_\nu \; \forall\nu \bigr\}.
\]
One may even fix any sequence of numbers $\vare_1,\ldots,\vare_l$, where $\vare_\nu=\pm 1$ for all $\nu$, and replace $\fD^+_{q^\star}$ or $\fD^-_{q^\star}$ with the domain
\[
\fD^\vare_{q^\star} = \bigl\{ (u,\varp,p,q) \bigm| p_\nu\in\fI_{\vare_\nu}\bmod 2\pi \; \forall\nu, \;\; q_\nu\neq q^\star_\nu \; \forall\nu \bigr\},
\]
where $\fI_1 = (-\pi/2,\pi/2)$ and $\fI_{-1} = (\pi/2,3\pi/2)$, so that $\vare_\nu\dot{q}_\nu\geq 0$ everywhere in $\fD^\vare_{q^\star}$, $1\leq\nu\leq l$.

If $l\geq 1$ and all the constants $\zeta_1,\ldots,\zeta_s$, $\xi_1,\ldots,\xi_l$, $\eta_1,\ldots,\eta_l$ are positive (so that $d=0$), then $\cT_{0,0,0}$ is the only Kronecker torus of \eqref{XourhatH} entirely contained in the domain
\[
\bigl\{ (u,\varp,p,q) \bigm| u_i\neq\pi \; \forall i, \;\; p_\nu\in(-\pi/2,\pi/2)\bmod 2\pi \; \forall\nu, \;\; q_\nu\neq\pi \; \forall\nu \bigr\}.
\]
So, in this case $\cT_{0,0,0}$ is strongly isolated.

\section{Reversible analogues}\label{reversible}

Both the Hamiltonian systems \eqref{XourH} and \eqref{XourhatH} are reversible with respect to the phase space involution
\[
\widetilde{G}: (u,\varp,p,q) \mapsto (u,-\varp,p,-q)
\]
of type $(s+2k+l,s+l)$, so that $\dim\Fix\widetilde{G}=s+l \geq 1$, $\codim\Fix\widetilde{G}=s+2k+l \geq \dim\Fix\widetilde{G}$, and $\codim\Fix\widetilde{G}-n=l < \dim\Fix\widetilde{G}$, where $n=s+2k$. However, $\widetilde{G}\bigl( \cT_{u^0,p^0,q^0} \bigr) = \cT_{u^0,p^0,-q^0}$, so that not all the $n$-tori \eqref{cT} are invariant under $\widetilde{G}$. In the case of the system \eqref{XourH}, a torus $\cT_{u^0,p^0,q^0}$ is invariant under $\widetilde{G}$ if and only if $q^0=0$. Consequently, the statement ``each torus $\cT_{u^0,p^0,q^0} \subset \fM$ is symmetric'' is valid if and only if all the numbers $\xi_1,\ldots,\xi_l$ are positive. In the case of the system \eqref{XourhatH}, a torus $\cT_{u^0,p^0,q^0}$ is invariant under $\widetilde{G}$ if and only if $q^0=-q^0$, i.e., if each component of $q^0$ is equal to either $0$ or $\pi$. Again, the statement ``each torus $\cT_{u^0,p^0,q^0} \subset \widehat{\fM}$ is symmetric'' holds if and only if all the numbers $\xi_1,\ldots,\xi_l$ are positive.

It is easy to construct a $G$-reversible counterpart of the system \eqref{XourH} for any non-negative integer values of $n$, $\dim\Fix G$, and $\codim\Fix G-n$, where $n$ is the dimension of symmetric Kronecker tori. Let $m$, $n$, $l$ be non-negative integers and consider the manifold \eqref{cK} equipped with the involution \eqref{inv} of type $(n+l,m)$. For any $u^0\in\mR^m$ and $q^0\in\mR^l$, consider the $n$-torus
\begin{equation}
\cT_{u^0,q^0} = \bigl\{ (u^0,\varp,q^0) \bigm| \varp\in\mT^n \bigr\}.
\label{TG}
\end{equation}
Since $G\bigl( \cT_{u^0,q^0} \bigr) = \cT_{u^0,-q^0}$, a torus $\cT_{u^0,q^0}$ is invariant under $G$ if and only if $q^0=0$.

Now let $\zeta_1,\ldots,\zeta_m$, $\xi_1,\ldots,\xi_l$ be arbitrary non-negative real constants and let $h:\mR^m\to\mR^n$ be an arbitrary smooth mapping. The system
\begin{equation}
\begin{aligned}
\dot{u}_i &= 0, \\
\dot{\varp}_\alpha &= h_\alpha(u), \\
\dot{q}_\nu &= \xi_\nu q_\nu^2 + \delta_{1\nu}l \sum_{i=1}^m \zeta_iu_i^2
\end{aligned}
\label{Xrev}
\end{equation}
(where $1\leq i\leq m$, $1\leq\alpha\leq n$, $1\leq\nu\leq l$) is reversible with respect to $G$. The term $l \sum_{i=1}^m \zeta_iu_i^2$ automatically vanishes for $l=0$. The key property of the system \eqref{Xrev} is that $\dot{q}_\nu\geq 0$ everywhere in the phase space $\cK$, $1\leq\nu\leq l$.

In the note \cite{S415}, we considered a similar system with $h(u)\equiv\omega\in\mR^n$, with $l\geq 1$, and with the equation for $\dot{q}_\nu$ of the form
\[
\dot{q}_\nu = \delta_{1\nu}\left( \sum_{\mu=1}^l q_\mu^2 + \sum_{i=1}^m u_i^2 \right),
\]
$1\leq\nu\leq l$. Our variables $m$ and $l$ play the roles of $\ell$ and $m+1$ in \cite{S415}, respectively.

All the conditionally periodic motions of the system \eqref{Xrev} fill up the manifold
\[
\fK = \bigl\{ (u,\varp,q) \bigm| l\zeta_iu_i=0 \; \forall i, \;\; \xi_\nu q_\nu=0 \; \forall\nu \bigr\}
\]
foliated by Kronecker $n$-tori of the form \eqref{TG}. Of course, always $\cT_{0,0}\subset\fK$, and the torus $\cT_{0,0}$ is symmetric. The frequency vector of a torus $\cT_{u^0,q^0} \subset \fK$ is $h(u^0)$, and this vector can be made equal to any prescribed vector $\omega\in\mR^n$ just by setting $h(u)\equiv\omega$. If $l=0$ then $\fK=\cK$. The system \eqref{Xrev} admits no conditionally periodic motions outside $\fK$.

These features of $\fK$ can be verified in exactly the same way as in Section~\ref{system}. If $(u,\varp,q)\in\fK$ then $\dot{u}=0$, $\dot{\varp}=h(u)$, $\dot{q}=0$. On the other hand, since $\dot{q}_\nu\geq 0$ everywhere in $\cK$, the recurrence property of conditionally periodic motions implies that $\dot{q}_\nu\equiv 0$ on Kronecker tori, $1\leq\nu\leq l$. Consequently, a point $(u,\varp,q)\notin\fK$ does not belong to any Kronecker torus of \eqref{Xrev} (symmetric or not and of any dimension) because $\dot{q}_1>0$ whenever $l\zeta_iu_i\neq 0$ for at least one $i$ and $\dot{q}_\nu>0$ whenever $\xi_\nu q_\nu\neq 0$, $1\leq\nu\leq l$.

The Kronecker $n$-tori $\cT_{u^0,q^0} \subset \fK$ constitute an analytic $d$-parameter family where $d = \dim\fK-n$. If $l=0$ then $d=m$. If $l\geq 1$ (i.e., if $\codim\Fix G>n$) then $d$ can take any integer value between $0$ and $m+l$; to be more precise, $d$ is the number of zero constants among $\zeta_i$, $\xi_\nu$ ($1\leq i\leq m$, $1\leq\nu\leq l$). The equality $d=0$ holds if and only if all the numbers $\zeta_i$, $\xi_\nu$ are positive in which case $\fK=\cT_{0,0}$, and $\cT_{0,0}$ is a unique Kronecker torus of the system \eqref{Xrev}. The equality $d=m+l$ occurs if and only if all the numbers $\zeta_i$, $\xi_\nu$ are equal to zero in which case $\fK=\cK$.

The symmetric Kronecker $n$-tori $\cT_{u^0,0}$ of the system \eqref{Xrev} are characterized by the condition $l\zeta_iu^0_i=0 \; \forall i$ and constitute an analytic $d_\ast$-parameter family where $d_\ast$ is determined as follows. If $l=0$ then $d_\ast=m$, and all the Kronecker $n$-tori constituting $\fK=\cK$ are symmetric. If $l\geq 1$ then $d_\ast$ is the number of zero constants among $\zeta_i$ ($1\leq i\leq m$) and can therefore take any integer value between $0$ and $m$. In all the cases, $d-d_\ast\leq l$.

\section{Compactified reversible analogues}\label{comprev}

The system \eqref{Xrev} can be compactified in the same way as the system \eqref{XourH}, cf.\ \cite{S415}. Consider the manifold \eqref{hatcK} equipped with the involution $G$ given by the same formula \eqref{inv} and having the same type $(n+l,m)$. For any $u^0\in\mT^m$ and $q^0\in\mT^l$, consider the $n$-torus $\cT_{u^0,q^0}$ given by the same expression \eqref{TG}. Since $G\bigl( \cT_{u^0,q^0} \bigr) = \cT_{u^0,-q^0}$, a torus $\cT_{u^0,q^0}$ is invariant under $G$ if and only if $q^0=-q^0$, i.e., if each component of $q^0$ is equal to either $0$ or $\pi$.

Now let again $\zeta_1,\ldots,\zeta_m$, $\xi_1,\ldots,\xi_l$ be arbitrary non-negative real constants and let $h:\mT^m\to\mR^n$ be an arbitrary smooth mapping. The system
\begin{equation}
\begin{aligned}
\dot{u}_i &= 0, \\
\dot{\varp}_\alpha &= h_\alpha(u), \\
\dot{q}_\nu &= \xi_\nu\tilde{q}_\nu^2 + \delta_{1\nu}l \sum_{i=1}^m \zeta_i\tilde{u}_i^2
\end{aligned}
\label{Xrevhat}
\end{equation}
(where $1\leq i\leq m$, $1\leq\alpha\leq n$, $1\leq\nu\leq l$, and the notation $\tilde{z}=\sin z$ is used) is reversible with respect to $G$, and $\dot{q}_\nu\geq 0$ everywhere in the phase space $\widehat{\cK}$, $1\leq\nu\leq l$.

The manifold
\[
\widehat{\fK} = \bigl\{ (u,\varp,q) \bigm| l\zeta_i\tilde{u}_i=0 \; \forall i, \;\; \xi_\nu\tilde{q}_\nu=0 \; \forall\nu \bigr\}
\]
is again foliated by Kronecker $n$-tori of the form \eqref{TG} (with $u^0\in\mT^m$ and $q^0\in\mT^l$), $\cT_{0,0}\subset\widehat{\fK}$ in all the cases, and the torus $\cT_{0,0}$ is symmetric. The frequency vector of a torus $\cT_{u^0,q^0} \subset \widehat{\fK}$ is $h(u^0)$, and this vector can be made equal to any prescribed vector $\omega\in\mR^n$ just by setting $h(u)\equiv\omega$. The dimension $n+d$ of the manifold $\widehat{\fK}$ is determined in exactly the same way as that of the manifold $\fK$ in Section~\ref{reversible}. In particular, if $l=0$ then $\widehat{\fK}=\widehat{\cK}$.

If $l\geq 1$ then $\sum_{\nu=1}^l \dot{q}_\nu > 0$ everywhere in $\widehat{\cK} \setminus \widehat{\fK}$. Like in Section~\ref{compact} and in contrast to the case of the system \eqref{Xrev}, it is, generally speaking, \emph{not} true that the system \eqref{Xrevhat} for $l\geq 1$ possesses no conditionally periodic motions outside $\widehat{\fK}$. Indeed, similarly to the example in Section~\ref{compact}, suppose that $m\geq 1$, $\sum_{i=1}^m \zeta_i > 0$ and choose an arbitrary point $u^0\in\mT^m$ such that $\chi = \sum_{i=1}^m \zeta_i\sin^2u^0_i > 0$. Consider the $(n+1)$-torus
\[
\bigl\{ (u^0,\varp,q) \bigm| q_2=\cdots=q_l=0 \bigr\} \not\subset \widehat{\fK}.
\]
This is a symmetric Kronecker torus of the system \eqref{Xrevhat} with the frequencies $h_1(u^0),\ldots,h_n(u^0),\varpi$, where $\varpi = \bigl[ l\chi(l\chi+\xi_1) \bigr]^{1/2}$.

Nevertheless, for any fixed $q^\star\in\mT^l$, no point $(u,\varp,q)\notin\widehat{\fK}$ belongs to a Kronecker torus of \eqref{Xrevhat} (symmetric or not and of any dimension) entirely contained in the domain
\[
\bigl\{ (u,\varp,q) \bigm| q_\nu\neq q^\star_\nu \; \forall\nu \bigr\}.
\]
This may be verified in exactly the same way as in Section~\ref{compact}.

If $l\geq 1$ and all the constants $\zeta_1,\ldots,\zeta_m$, $\xi_1,\ldots,\xi_l$ are positive (so that $d=0$), then $\cT_{0,0}$ is the only Kronecker torus of \eqref{Xrevhat} entirely contained in the domain
\[
\bigl\{ (u,\varp,q) \bigm| u_i\neq\pi \; \forall i, \;\; q_\nu\neq\pi \; \forall\nu \bigr\}.
\]
So, in this case $\cT_{0,0}$ is strongly isolated.

The symmetric Kronecker $n$-tori $\cT_{u^0,q^0}$ of the system \eqref{Xrevhat} make up an $(n+d_\ast)$-dimensional submanifold of the manifold $\widehat{\fK}$, where $d_\ast$ is determined in exactly the same way as in Section~\ref{reversible}.

\section*{Declaration of interest}

Declarations of interest: none.

\section*{Acknowledgments}

I am grateful to B.~Fayad for fruitful correspondence and sending me the breakthrough preprint \cite{FF01575} prior to submission to arXiv.


\begin{thebibliography}{99}

\bibitem{A1969}
J.F.~Adams, Lectures on Lie Groups, W.A.~Benjamin, Inc., New York, 1969.

\bibitem{AKN2006}
V.I.~Arnold, V.V.~Kozlov, A.I.~Neishtadt, Mathematical Aspects of Classical and Celestial Mechanics, third ed. (Encyclopaedia of Mathematical Sciences, vol.~3), Springer, Berlin, 2006.

\bibitem{B2005}
R.~Baer, Linear Algebra and Projective Geometry, second ed., Dover Publ., Mineola, NY, 2005.

\bibitem{BDP223}
E.M.~Barbaresco, P.E.~Desideri, P.L.Q.~Pergher, Involutions whose fixed set has three or four components:\ a small codimension phenomenon, Math. Scand. 110 (2) (2012) 223--234.

\bibitem{BZ365}
M.S.~Borman, F.~Zapolsky, Quasimorphisms on contactomorphism groups and contact rigidity, Geom. Topol. 19 (1) (2015) 365--411.

\bibitem{B42}
A.~Bounemoura, Non-degenerate Liouville tori are KAM stable, Adv. Math. 292 (2016) 42--51.

\bibitem{B21}
A.~Bounemoura, Some instability properties of resonant invariant tori in Hamiltonian systems, Math. Res. Lett. 24 (1) (2017) 21--35.

\bibitem{BHN355}
H.W.~Broer, J.~Hoo, V.~Naudot, Normal linear stability of quasi-periodic tori, J. Differential Equations 232 (2) (2007) 355--418.

\bibitem{BH191}
H.W.~Broer, G.B.~Huitema, Unfoldings of quasi-periodic tori in reversible systems, J. Dynam. Differential Equations 7 (1) (1995) 191--212.

\bibitem{BHS1996}
H.W.~Broer, G.B.~Huitema, M.B.~Sevryuk, Quasi-Periodic Motions in Families of Dynamical Systems. Order amidst Chaos (Lecture Notes in Mathematics, vol.~1645), Springer, Berlin, 1996.

\bibitem{BHTB1990}
H.W.~Broer, G.B.~Huitema, F.~Takens, B.L.J.~Braaksma, Unfoldings and Bifurcations of Quasi-Periodic Tori (Mem. Amer. Math. Soc., vol.~83, no.~421), Amer. Math. Soc., Providence, RI, 1990.

\bibitem{BS249}
H.W.~Broer, M.B.~Sevryuk, KAM theory:\ quasi-periodicity in dynamical systems, in: H.W.~Broer, B.~Hasselblatt, F.~Takens (Eds.), Handbook of Dynamical Systems. Vol.~3, Elsevier B.V., Amsterdam, 2010, pp.~249--344.

\bibitem{B1}
J.~Butterfield, On symplectic reduction in classical mechanics, in: J.~Butterfield, J.~Earman (Eds.), Handbook of the Philosophy of Science. Philosophy of Physics, Elsevier B.V., Amsterdam, 2007, pp.~1--131.

\bibitem{dlL175}
R.~de la Llave, A tutorial on KAM theory, in: A.~Katok, R.~de la Llave, Ya.~Pesin, H.~Weiss (Eds.), Smooth Ergodic Theory and Its Applications (Proc. Sympos. Pure Math., vol.~69), Amer. Math. Soc., Providence, RI, 2001, pp.~175--292.

\bibitem{DP3119}
P.E.~Desideri, P.L.Q.~Pergher, Involutions fixing many components:\ a small codimension phenomenon, J. Fixed Point Theory Appl. 19 (4) (2017) 3119--3126.

\bibitem{DK2000}
J.J.~Duistermaat, J.A.C.~Kolk, Lie Groups, Springer, Berlin, 2000.

\bibitem{D2014}
H.S.~Dumas, The KAM Story. A Friendly Introduction to the Content, History, and Significance of Classical Kolmogorov--Arnold--Moser Theory, World Sci. Publ., Hackensack, NJ, 2014.

\bibitem{EM2002}
Ya.M.~Eliashberg, N.M.~Mishachev, Introduction to the $h$-Principle, Amer. Math. Soc., Providence, RI, 2002.

\bibitem{EFK1733}
L.H.~Eliasson, B.~Fayad, R.~Krikorian, Around the stability of KAM tori, Duke Math. J. 164 (9) (2015) 1733--1775.

\bibitem{FF01575}
G.~Farr\'e, B.~Fayad, Instabilities for analytic quasi-periodic invariant tori, 2019, https://arxiv.org/abs/1912.01575.

\bibitem{F93}
F.~Fass\`o, Superintegrable Hamiltonian systems:\ geometry and perturbations, Acta Appl. Math. 87 (1--3) (2005) 93--121.

\bibitem{F09059}
B.~Fayad, Lyapunov unstable elliptic equilibria, second version, 2019, https://arxiv.org/abs/1809.09059.

\bibitem{FK1905}
B.~Fayad, R.~Krikorian, Some questions around quasi-periodic dynamics, in: B.~Sirakov, P.~Ney de Souza, M.~Viana (Eds.), Proceedings of the International Congress of Mathematicians (Rio de Janeiro, 2018), vol.~III. Invited lectures, World Sci. Publ., Hackensack, NJ, 2018, pp.~1905--1930.

\bibitem{FS67}
B.~Fayad, M.~Saprykina, Isolated elliptic fixed points for smooth Hamiltonians, in: A.~Katok, Ya.~Pesin, F.~Rodriguez Hertz (Eds.), Modern Theory of Dynamical Systems. A Tribute to Dmitry Victorovich Anosov (Contemp. Math., vol.~692), Amer. Math. Soc., Providence, RI, 2017, pp.~67--82.

\bibitem{F1521}
J.~F\'ejoz, D\'emonstration du `th\'eor\`eme d'Arnold' sur la stabilit\'e du syst\`eme plan\'etaire (d'apr\`es Herman), Ergodic Theory Dynam. Systems 24 (5) (2004) 1521--1582.

\bibitem{HLW2006}
E.~Hairer, Ch.~Lubich, G.~Wanner, Geometric Numerical Integration. Structure-Preserving Algorithms for Ordinary Differential Equations, second ed., Springer, Berlin, 2006.

\bibitem{H79}
H.~Han{\ss}mann, Perturbations of superintegrable systems, Acta Appl. Math. 137 (2015) 79--95.

\bibitem{H989}
M.-R.~Herman, Exemples de flots hamiltoniens dont aucune perturbation en topologie $C^\infty$ n'a d'orbites p\'eriodiques sur un ouvert de surfaces d'\'energies, C. R. Acad. Sci. Paris S\'er.~I Math. 312 (13) (1991) 989--994.

\bibitem{H49}
M.-R.~Herman, Diff\'erentiabilit\'e optimale et contre-exemples \`a la fermeture en topologie $C^\infty$ des orbites r\'ecurrentes de flots hamiltoniens, C. R. Acad. Sci. Paris S\'er.~I Math. 313 (1) (1991) 49--51.

\bibitem{H797}
M.~Herman, Some open problems in dynamical systems, in: Proceedings of the International Congress of Mathematicians (Berlin, 1998), vol.~II. Sections 1--9, Doc. Math. 1998, Extra vol.~II, pp.~797--808.

\bibitem{KP2003}
T.~Kappeler, J.~P\"oschel, KdV \& KAM, Springer, Berlin, 2003.

\bibitem{K2018}
Khanickus (https://mathoverflow.net/users/85369/khanickus), An answer to:\ Isolated periodic trajectories of Hamiltonian systems, 2018, https://mathoverflow.net/q/289725.

\bibitem{K259}
S.B.~Kuksin, An infinitesimal Liouville--Arnold theorem as a criterion of reducibility for variational Hamiltonian equations, Chaos Solitons Fractals 2 (3) (1992) 259--269.

\bibitem{KS811}
A.V.~Kurov, G.A.~Sardanashvili, Globally superintegrable Hamiltonian systems, Theoret. Math. Phys. 191 (3) (2017) 811--826.

\bibitem{MP643}
T.~Mitev, G.~Popov, Gevrey normal form and effective stability of Lagrangian tori, Discrete Contin. Dyn. Syst. Ser.~S 3 (4) (2010) 643--666.

\bibitem{O1993}
P.J.~Olver, Applications of Lie Groups to Differential Equations, second ed., Springer, New York, 1993.

\bibitem{P380}
I.O.~Parasyuk, On the persistence of multidimensional invariant tori of Hamiltonian systems, Ukrain. Math. J. 36 (4) (1984) 380--385.

\bibitem{PF280}
P.L.Q.~Pergher, F.G.~Figueira, Dimensions of fixed point sets of involutions, Arch. Math. (Basel) 87 (3) (2006) 280--288.

\bibitem{P707}
J.~P\"oschel, A lecture on the classical KAM theorem, in: A.~Katok, R.~de la Llave, Ya.~Pesin, H.~Weiss (Eds.), Smooth Ergodic Theory and Its Applications (Proc. Sympos. Pure Math., vol.~69), Amer. Math. Soc., Providence, RI, 2001, pp.~707--732.

\bibitem{QS757}
G.R.W.~Quispel, M.B.~Sevryuk, KAM theorems for the product of two involutions of different types, Chaos 3 (4) (1993) 757--769.

\bibitem{S2007}
M.R.~Sepanski, Compact Lie Groups, Springer, New York, 2007.

\bibitem{S351}
M.B.~Sevryuk, KAM-stable Hamiltonians, J. Dynam. Control Systems 1 (3) (1995) 351--366.

\bibitem{S1113}
M.B.~Sevryuk, The classical KAM theory at the dawn of the twenty-first century, Mosc. Math. J. 3 (3) (2003) 1113--1144.

\bibitem{St177}
M.B.~Sevryuk, KAM tori:\ persistence and smoothness, Nonlinearity 21 (10) (2008) T177--T185.

\bibitem{S435}
M.B.~Sevryuk, KAM theory for lower dimensional tori within the reversible context~2, Mosc. Math. J. 12 (2) (2012) 435--455.

\bibitem{S137}
M.B.~Sevryuk, Quasi-periodic perturbations within the reversible context~2 in KAM theory, Indag. Math. 23 (3) (2012) 137--150.

\bibitem{S599}
M.B.~Sevryuk, Whitney smooth families of invariant tori within the reversible context~2 of KAM theory, Regul. Chaotic Dyn. 21 (6) (2016) 599--620.

\bibitem{S603}
M.B.~Sevryuk, Families of invariant tori in KAM theory:\ interplay of integer characteristics, Regul. Chaotic Dyn. 22 (6) (2017) 603--615.

\bibitem{S415}
M.B.~Sevryuk, Integrable Hamiltonian systems with a periodic orbit or invariant torus unique in the whole phase space, Arnold Math. J. 4 (3--4) (2018) 415--422.

\bibitem{S6215}
C.~Sim\'o, Some questions looking for answers in dynamical systems, Discrete Contin. Dyn. Syst. 38 (12) (2018) 6215--6239.

\bibitem{T2012}
G.~Teschl, Ordinary Differential Equations and Dynamical Systems, Amer. Math. Soc., Providence, RI, 2012.

\bibitem{T12851}
F.~Trujillo, Lyapunov unstable elliptic equilibria for Hamiltonians with two degrees of freedom, 2019, https://arxiv.org/abs/1912.12851.

\bibitem{W329}
A.~Weinstein, Symplectic manifolds and their Lagrangian submanifolds, Adv. Math. 6 (3) (1971) 329--346.

\end{thebibliography}
\end{document}